\documentclass[12pt]{article}
\usepackage{amsfonts,amscd,a4}
\usepackage{color}
\usepackage{graphicx}
\usepackage{theorem}
\newtheorem{theo}{Theorem}
\newtheorem{prop}[theo]{Proposition}
\newtheorem{lemma}[theo]{Lemma}
\newtheorem{coro}[theo]{Corollary}

{\theorembodyfont{\rm}
\newtheorem{example}[theo]{Example}
}

\newcommand{\cA}{{\mathcal A}}
\newcommand{\cB}{{\mathcal B}}
\newcommand{\cC}{{\mathcal C}}

\newcommand{\cI}{{\mathcal I}}
\newcommand{\cJ}{{\mathcal J}}
\newcommand{\cK}{{\mathcal K}}

\newcommand{\sC}{{\mathbb C}}
\newcommand{\sD}{{\mathbb D}}

\newcommand{\sN}{{\mathbb N}}

\newcommand{\sR}{{\mathbb R}}

\newcommand{\sT}{{\mathbb T}}

\newcommand{\qed}{\rule{1ex}{1ex}}
\newcommand{\alg}{\mbox{\rm alg} \,}
\newcommand{\ann}{\mbox{\rm Ann} \,}

\newcommand{\id}{\mbox{\rm Id} \,}
\newcommand{\im}{\mbox{\rm im} \,}

\newcommand{\rad}{{\rm Rad} \,}

\begin{document}
\title{On Moore-Penrose ideals}
\author{Steffen Roch}
\date{}
\maketitle
\begin{abstract}
Every bounded linear operator on a Hilbert space which is invertible modulo compact operators has a closed range and is, thus, generalized invertible. We consider the analogue question in general $C^*$-algebras and describe the closed ideals (called Moore-Penrose ideals in what follows) with the property that whenever an element is invertible modulo that ideal, then it is generalized invertible. In particular, we will see that the class of Moore-Penrose ideals coincides with the class of the dual ideals. Finally, we study some questions related with the projection lifting property of Moore-Penrose ideals. 
\end{abstract}
{\bf 2010 AMS-MSC:} 46L05, 47A05
\section{Introduction}
The ideal $K(H)$ of the compact operators on a Hilbert space $H$ owns the following property: Whenever an operator $A \in L(H)$, the algebra of the bounded linear operators on $H$,  is invertible modulo $K(H)$, then its range $\im A = AH$ is closed in $H$. (Actually, $A$ is even a Fredholm operator, i.e. it has a finite-dimensional kernel $\ker A = A^{-1} (0)$ and a finite-dimensional cokernel $H/\im A$.)

The closedness of the range of $A$, also termed as the normal solvability of $A$, is equivalent to the generalized invertibility of $A$, i.e. to the existence of an operator $B \in L(H)$ such that $ABA = A$. In general, the generalized inverse $B$ of $A$ is not unique; but if $A$ has a generalized inverse, then it also has a generalized inverse $C$ with the following properties: $ACA = A$, $CAC = C$, and $AC$ and $CA$ are self-adjoint. A generalized inverse with these properties is unique; it is called the Moore-Penrose inverse of $A$ and denoted by $A^\dagger$. Thus, the terms {\em generalized invertibility} and {\em Moore-Penrose invertibility} (MPI for brevity) are synonymous.

These facts hold more general when $L(H)$ is replaced by an arbitrary $C^*$-algebra, in which setting generalized or Moore-Penrose invertibility often serves as a substitute of the closed range property of operators. The aim of this paper is to characterize closed ideals $\cJ$ of a $C^*$-algebra $\cA$ with the following property, mimicking the property of the ideal of the compact operators observed above: {\em whenever $a \in \cA$ is invertible modulo $\cJ$, then $a$ is Moore-Penrose invertible}. In this case we call $\cJ$ a {\em Moore-Penrose ideal} (MP-ideal for short) of $\cA$.

To state our results we have to introduce some more notations. A $C^*$-algebra $\cA$ is called {\em elementary} if it is isomorphic to the algebra of the compact operators on some Hilbert space, and a $C^*$-algebra is {\em dual} if it is isomorphic to a direct sum of elementary algebras. If $\cJ$ a closed ideal of $\cA$ which is elementary (respective dual) when considered as a $C^*$-algebra, then we call $\cJ$ an {\em elementary} (respective a {\em dual}) ideal \index{ideal!dual} of $\cA$. It is easy to see that every dual ideal $\cJ$ is generated (as a $C^*$-algebra) by its elementary ideals, $\cK_t$ with $t \in T$, say. Since every closed ideal of $\cJ$ is also a closed ideal of $\cA$, $\cJ$ can be identified with the smallest closed ideal of $\cA$ which contains all elementary ideals ideals $\cK_t$. See \cite{Ber1} for an overview on dual algebras.

Further we say that an ideal $\cJ$ of a unital $C^*$-algebra has the {\em compact spectral property} if every self-adjoint element of $\cJ$ has 0 as the only possible accumulation point of its spectrum. Clearly, this notion is motivated by the fact that every compact operator on a Hilbert space has $0$ as only possible accumulation point of its spectrum.

Finally, we call a non-zero element $k$ of a $C^*$-algebra $\cA$ an {\em element of algebraic rank one} if, for every $a \in \cA$, there is a complex number $\alpha$ such that $kak = \alpha k$. The smallest closed ideal of $\cA$ which contains all elements of algebraic rank one is denoted by $\cC(\cA)$. Its elements are called {\em algebraically compact}.
\begin{theo} \label{tb3.90}
The following assertions are equivalent for a proper closed ideal $\cJ$ of a unital $C^*$-algebra $\cB \!:$ \\[1mm]
$(a)$ $\cJ$ is a dual ideal. \\[1mm]
$(b)$ $\cJ$ is a Moore-Penrose ideal. \\[1mm]
$(c)$ $\cJ$ has the compact spectral property. \\[1mm]
$(d)$ $\cJ$ is $^*$-isomorphic to a $C^*$-subalgebra of $K(H)$ for a Hilbert space $H$. \\[1mm]
$(e)$ $\cC(\cJ) = \cJ$.
\end{theo}
The implications $(a) \Rightarrow (b)$, $(b) \Rightarrow (c)$ and $(c) \Rightarrow (e)$ are shown in Sections \ref{ssb1.4}, \ref{ssb1.5} and \ref{ssb1.6}, respectively. The remaining implications are subject to Section \ref{ssb1.6b}. Of course, not all of the implications in Theorem \ref{tb3.90} are new; for example, the implication $(d) \Rightarrow (a)$ is Theorem 1.4.5 in \cite{Arv1}.

The following theorem summarizes some kown facts on Moore-Penrose invertibility in $C^*$-algebras for later reference.
\begin{theo} \label{t31.2}
Let $\cB$ be a $C^*$-algebra with identity $e$. The following conditions are equivalent for an element $a$ of $\cB$: \\[1mm]
$(a)$ $a$ is generalized invertible. \\[1mm]
$(b)$ $a$ is Moore-Penrose invertible. \\[1mm]
$(c)$ $a^*a$ is invertible, or $0$ is an isolated point of the spectrum $\sigma(a^*a)$ of $a^*a$. \\[1mm]
$(d)$ There is a projection $p \in \alg (e, \, a^*a)$ such that $ap = 0$ and $a^*a + p$ is invertible. (Here, $\alg (e, \, a^*a)$ the smallest closed subalgebra of $\cB$ which contains $e$ and $a^*a$, and $p$ is called a projection if $p$ is self-adjoint and $p^2 = p$.) \\[1mm]
$(e)$ There is a projection $q \in \cB$ such that $aq = 0$ and $a^*a + q$ is invertible. \\[1mm]
If one of these conditions is satisfied, then the projection $q$ is uniquely determined, $q = e - a^\dagger a$, and $a^\dagger = (a^*a + q)^{-1} a^*$. Moreover,
\[
\|a^\dagger\|^2 = \sup \, \{ 1/\lambda : \lambda \in \sigma (a^*a) \setminus \{0\} \}.
\]
$(f)$ If $\cC$ is a unital $C^*$-subalgebra of $\cB$ and $c \in \cC$ is Moore-Penrose invertible in $\cB$, then $c^\dagger \in \cC$ (inverse closedness).
\end{theo}
\section{Dual ideals and lifting theorems} \label{ssb1.1}
%
%
%
Our study of dual ideals is dominated by a lifting theorem. Lifting theorems can be used to reduce invertibility problems in a $C^*$-algebra $\cA$ to invertibility problems in a suitable quotient algebra $\cA/\cJ$ (in some sense, they ``measure`` the difference between invertibility in $\cA$ and in $\cA/\cJ$). They hold in quite general contexts (see Section 6.3 in \cite{RSS2} for purely algebraic and Banach algebraic versions of that theorem), but they work particularly well when $\cJ$ is a dual algebra.

Let $\cA$ and $\cB$ be unital $C^*$-algebras, $\cJ$ a closed ideal of $\cA$, and $W : \cA \to \cB$ a unital $^*$-homomorphism. We say that $W$ {\em lifts the ideal} $\cJ$ if the  restriction of $W$ to $\cJ$ is injective. Note that for every closed ideal $\cJ$ there is a canonical homomorphism which lifts $\cJ$. Indeed, the annulator
\[
\ann \cJ := \{a \in \cA : a \cJ \cup \cJ a \subseteq \rad \cA \}
\]
of $\cJ$ is a closed ideal of $\cA$, and the canonical homomorphism from $\cA$ to $\cA/\ann \cJ$ lifts the ideal $\cJ$.

Here is a general version of the lifting theorem for $C^*$-algebras (Theorem 6.3.6 in \cite{RSS2}).
\begin{theo}[Lifting theorem] \label{ta6.8}
Let $\cA$ be a $C^*$-algebra with identity. For every element $t$ of a certain set $T$, let $\cJ_t$ be a closed ideal of $\cA$ which is lifted by a unital $^*$-homomorphism $W_t$ from $\cA$ into a unital $^*$-algebra $\cB_t$. Let further $\cJ$ stand for the smallest closed ideal of $\cA$ which contains all ideals $\cJ_t$. Then an element $a \in \cA$ is invertible if and only if the coset $a + \cJ$ is invertible in $\cA/\cJ$ and if all elements $W_t(a)$ are invertible in $\cB_t$.
\end{theo}
The family $(W_t)_{t \in T}$ of homomorphisms induces a product homomorphism $W$ from $\cA$ into the product of the family $(\cB_t)_{t \in T}$ via
\begin{equation} \label{ea6.4}
W : a \mapsto (t \mapsto W_t(a)).
\end{equation}
Thus, the lifting theorem states that the homomorphisms $W_t$ and the ideals $\cJ_t$ can be glued to a homomorphism $W$ and to an ideal $\cJ$, respectively, such that $W$ lifts $\cJ$.

Now we return to the context of dual ideals. Our next goal is to show that every dual ideal of a $C^*$-algebra can be lifted and to specify Theorem \ref{ta6.8} to this context. We start with lifting one elementary ideal.
\begin{prop} \label{pb1.5}
Let $\cA$ be a unital $C^*$-algebra and $\cJ$ an elementary ideal of $\cA$. Then there is an irreducible representation of $\cA$ which lifts $\cJ$. This representation is unique up to unitary equivalence, and it maps the elements of $\cJ$ to compact operators.
\end{prop}
{\bf Proof.} Let $W$ be a $^*$-isomorphism from $\cJ$ to $K(H)$. Then $W$ is the only irreducible representation of $\cJ$ up to unitary equivalence, and there is a unique (again, up to unitary equivalence) extension of $W$ to an irreducible representation of $\cA$. This extension lifts $\cJ$ by its construction, and it is unital since every irreducible representation of unital algebras is unital. \hfill \qed
\begin{theo}[Lifting theorem for dual ideals] \label{tb1.6}
Let $\cA$ be a unital $C^*$-al\-gebra. For every element $t$ of a set $T$, let $\cJ_t$ be an elementary ideal of $\cA$ such that $\cJ_s \neq \cJ_t$ whenever $s \neq t$, and let $W_t : \cA \to L(H_t)$ denote the irreducible representation of $\cA$ which lifts $\cJ_t$ $($which exists and is unique by Proposition $\ref{pb1.5})$. Let finally $\cJ$ stand for the smallest closed ideal of $\cA$ which contains all ideals $\cJ_t$. Then $\cJ$ is a dual ideal of $\cA$, and the assertion of the lifting Theorem $\ref{ta6.8}$ can be completed as follows: \\[1mm]
$(a)$ If $t_1, \, \ldots, \,  t_m \in T$ and $t_i \neq t_j$ for $i \neq j$, then $(\cJ_{t_1} + \ldots + \cJ_{t_{m-1}}) \cap \cJ_{t_m} = \{0\}$. \\[1mm]
$(b)$ The separation property \index{separation property} holds, i.e. $W_s (\cJ_t) = \{0\}$ whenever $s \neq t$. \\[1mm]
$(c)$ If $j \in \cJ$, then the operator $W_t(j)$ is compact for every $t \in T$, and the set of all $t \in T$ with $\|W_t(j)\| \ge \varepsilon$ is finite for every $\varepsilon > 0$. \\[1mm]
$(d)$ If $a \in \cA$ and the coset $a + \cJ$ is invertible, then the operator $W_t(a)$ is Fredholm, and the set of all $t \in T$ for which $W_t(a)$ is is not invertible is finite. \\[1mm]
$(e)$ The product $W := \Pi_{t \in T} W_t : \cA \to \Pi_{t \in T} L(H_t)$ is a $^*$-isomorphism between $\cJ$ and the direct sum $\oplus_{t \in T} K(H_t)$.
\end{theo}
{\bf Proof.} $(a)$ First we show that $\cJ_s \cap \cJ_t = \{0\}$ whenever $s \neq t$. Clearly, $\cJ_s \cap \cJ_t$ is a closed ideal of $\cJ_s$. Since $\cJ_s$ is an elementary algebra, its only closed ideals are $\{0\}$ and $\cJ_s$. If $\cJ_s \cap \cJ_t = \cJ_s$, then $\cJ_s \subseteq \cJ_t$. Since $\cJ_s \neq 0$ and $\cJ_t$ is elementary, this implies $\cJ_s = \cJ_t$, a contradiction. Thus, $\cJ_s \cap \cJ_t = \{0\}$. This equality also implies that $\cJ$ is a dual algebra.

Now consider the case when $m \ge 2$. Clearly, $(\cJ_{t_1} + \ldots + \cJ_{t_{m-1}}) \cap \cJ_{t_m}$ is a closed ideal of $\cJ_{t_m}$. Since $\cJ_{t_m}$ is an elementary algebra, its only closed ideals are $\{0\}$ and $\cJ_{t_m}$. Thus, assuming that the assertion is wrong, we get $\cJ_{t_m} \subseteq \cJ_{t_1} + \ldots + \cJ_{t_{m-1}}$. Let $p \in \cJ_{t_m}$ be a non-zero projection (which exists since this ideal isomorphic to $K(H)$ for a non-zero Hilbert space $H$), and write $p = k_{t_1} + \ldots + k_{t_{m-1}}$ with elements $k_{t_i} \in \cJ_{t_i}$. Multiplying this identity by $p$ yields $p = p k_{t_1} + \ldots + p k_{t_{m-1}}$. Since $p \neq 0$, there is an $i$ such that $p k_{t_i} \neq 0$. Thus, $p k_{t_i}$ is a non-zero element of $\cJ_{t_m} \cap \cJ_{t_i}$, which is impossible by what we have already shown. \\[1mm]
$(b)$ Let $s, \, t \in T$ with $s \neq t$. Since $\cJ_s \cap \cJ_t = \{0\}$ by part $(a)$, every element $j \in \cJ_s + \cJ_t$ has a unique representation as $j = j_s + j_t$ with $j_s \in \cJ_s$ and $j_t \in \cJ_t$. This implies that the mapping
\[
\hat{W}_s : \cJ_s + \cJ_t \to K(H_s), \quad j = j_s + j_t \mapsto W_s(j_s)
\]
is correctly defined. This mapping is an irreducible representation of $\cJ_s + \cJ_t$ which coincides with $W_s$ on $\cJ_s$. Furthermore, $\hat{W}_s (\cJ_t) = \{0\}$.

Since the irreducible extension of $W_s$ from $\cJ_s$ onto $\cA$ is unique up to unitary equivalence, $\hat{W}_s$ is unitarily equivalent to the restriction to the ideal $\cJ_s + \cJ_t$ of any irreducible extension of $W_s$. Since equivalent representations have the same kernels, and since $\cJ_t$ lies in the kernel of $\hat{W}_s$, we get the assertion. \\[1mm]
$(c)$ Let $j \in \cJ$ and $\varepsilon > 0$. By the definition of $\cJ$, there exist a finite subset $\{t_1, \, \ldots, \, t_n\}$ of $T$ and elements $j_{t_i} \in \cJ_{t_i}$ such that
\begin{equation} \label{eb1.7}
j = j_{t_1} + \ldots + j_{t_n} + j^\prime \quad \mbox{with} \; j^\prime \in \cJ \; \mbox{and} \; \|j^\prime\| < \varepsilon.
\end{equation}
Let $t \in T$. Applying $W_t$ to both sides of (\ref{eb1.7}) and taking into account the separation property $(b)$, we get the existence of a compact operator $J_{t, \varepsilon} \in K(H_t)$ such that
\[
W_t(j) = J_{t, \varepsilon} + W_t(j^\prime) \quad \mbox{with} \quad \|W_t(j^\prime)\| < \varepsilon.
\]
The compactness of $W_t(j)$ follows since $K(H_t)$ is closed. Moreover, if $t \in T \setminus \{t_1, \, \ldots, \, t_n\}$, the same reasoning shows that $\|W_t(j)\| = \|W_t(j^\prime)\| < \varepsilon$. \\[1mm]
$(d)$ If the coset $a + \cJ$ is invertible, then there are elements $b \in \cA$ and $j, \, k \in \cJ$ such that $ab = e + j$ and $ba = e + k$. Applying the homomorphism $W_t$ to both equalities and taking into account assertion $(c)$ we obtain the Fredholm property of $W_t(a)$ for every $t \in T$. For the second assertion, write $j$ as in (\ref{eb1.7}) with some $\varepsilon < 1$. Applying $W_t$ for $t \not\in \{t_1, \, \ldots, \, t_n\}$ to the equation $ab = e + j$ one gets $W_t(a) W_t(b) = I_t + W_t(j^\prime)$ with $\|W_t(j^\prime)\| < 1$. Similarly, $W_t(b) W_t(a) = I_t + W_t(k^\prime)$ with $\|W_t(k^\prime)\| < 1$. Hence, $W_t(a)$ is invertible by Neumann series for all but finitely many $t$. \\[1mm]
$(e)$ The product homomorphism $W$ maps the ideal $\cJ$ to the direct product of the ideals $K(H_t)$ as we have already seen, and $W$ is injective due to Corollary 6.3.3 in \cite{RSS2}. Hence, $\cJ$ is $^*$-isomorphic to $W(\cJ)$. We show that the image of $\cJ$ under $W$ is just the direct sum of the ideals $K(H_t)$ and, hence, a dual algebra. Let $r \in W(\cJ)$, and let $j$ denote the unique element of $\cJ$ with $W(j) = r$. Given $\varepsilon > 0$, there is a decomposition of $j$ as in (\ref{eb1.7}). This decomposition immediately shows that $k$ belongs to the sum of the $K(H_t)$. On the other hand, the separation property entails that $W(\cJ_t) = K(H_t)$; hence, $W(\cJ)$ cannot be smaller than this sum. \hfill \qed \\[3mm]
If $\cJ$ is $^*$-isomorphic to the sum of its elementary ideals $\cJ_t$ and $W_t : \cJ_t \to K(H_t)$ is the (unique up to unitary equivalence) irreducible representation of $\cJ_t$, then $\cJ_s \cap \cJ_t = \{0\}$ if $s \neq t$, and $\cJ$ is the smallest closed ideal of $\cA$ which contains all ideals $\cJ_t$. It is also clear that every mapping $W_t$ has a unique (up to unitary equivalence) extension to an irreducible representation of $\cA$ on $H_t$, which lifts the ideal $\cJ_t$. Thus, it is exactly the dual ideals which satisfy the assumptions of Theorem \ref{tb1.6}.

Our final goal in this section is spectral and strictly spectral families of lifting homomorphisms. Given a unital $C^*$-algebra $\cA$ and a family $\{W_t\}_{t \in T}$ of unital homomorphisms $W_t$ from $\cA$ into unital $C^*$-algebras $\cB_t$, we call this family {\em spectral} if an element $a \in \cA$ is invertible if and only if all elements $W_t(a)$ are invertible and if the norms of their inverses are uniformly bounded with respect to $t \in T$. The family $\{W_t\}_{t \in T}$ is called {\em strictly spectral} if $a \in \cA$ is invertible if and only if every $W_t(a)$ is invertible. A family $\{W_t\}_{t \in T}$ is spectral if and only if the equality 
\begin{equation} \label{e100}
\|a\| = \sup_{t \in T} \|W_t(a)\|
\end{equation}
holds for every $a \in \cA$; it is strictly spectral if the maximum is attained in (\ref{e100}) for every $a \in \cA$. The first assertion is immediate from the fact that every injective $^*$-homomorphism is isometric; for the second one see, e.g., Theorem 2.2.11 in \cite{RSS2} (where strictly spectral families were called {\em sufficient families}).

The following theorem states converses of assertions $(c)$ and $(d)$ of the Lifting theorem \ref{tb1.6}. This result is then employed to characterize the ideals which can be lifted by a spectral family of homomorphisms.
\begin{theo} \label{tb1.8}
Let the notation and hypotheses be as in Theorem $\ref{tb1.6}$. Let $a, \, j \in \cA$. \\[1mm]
$(a)$ If the family $(W_t)_{t \in T}$ is spectral for $\cA$, then $j \in \cJ$ if and only if all operators $W_t(j)$ are compact and, for each $\varepsilon > 0$, the number of $t \in T$ with $\|W_t(j)\| \ge \varepsilon$ is finite. \\[1mm]
$(b)$ If the family $(W_t)_{t \in T}$ is spectral for $\cA$, then the coset $a + \cJ$ is invertible in $\cA/\cJ$ if and only if all operators $W_t(a)$ are Fredholm, if all but finitely many of them are invertible, and if the norms $\|(W_t(a) + K(H_t))^{-1}\|$ are uniformly bounded with respect to $t \in T$. \\[1mm]
$(c)$ If the family $(W_t)_{t \in T}$ is strictly spectral for $\cA$, then the coset $a + \cJ$ is invertible in $\cA/\cJ$ if and only if all operators $W_t(a)$ are Fredholm and if all but finitely many of them are invertible.
\end{theo}
{\bf Proof.} $(a)$ Due to Theorem \ref{tb1.6} $(c)$, we only have to show the "if"-part of the assertion. For $n \in \sN$, let $T_n := \{t \in T : \|W_t(j)\| \ge 1/n\}$. The sets $T_n$ are finite by hypothesis. For each $t \in T_n$, there is a (uniquely determined) element $k_t \in \cJ_t$ such that $W_t(k_t) = W_t(j)$. Set $j_n := \sum_{t \in T_n} k_t$. Then $\|W_t(j - j_n)\| < 1/n$ for every $t \in T$. Since $(W_t)_{t \in T}$ is a spectral family, we conclude via (\ref{e100}) that $\|j - j_n\| \le 1/n$. Since $j_n \in \cJ$ and $\cJ$ is closed, this implies that $j \in \cJ$. \\[1mm]
$(b)$ From Theorem \ref{tb1.6} $(d)$ we infer that if $a + \cJ$ is invertible then all operators $W_t(a)$ are Fredholm and only a finite number of them are not invertible. For the second assertion, write $ab = e + j$ with certain elements $b \in \cA$ and $j \in \cJ$. This identity implies $(W_t(a) + K(H_t))^{-1} = W_t(b) + K(H_t)$ and the estimate $\|W_t(b) + K(H_t)\| \le \|b\|$ for every $t \in T$. \\[1mm]
For the reverse implication, let all operators $W_t(a)$ be Fredholm and let $S \subseteq T$ be a finite set such that $W_t(a)$ be invertible for $t \in T \setminus S$. Then all operators $W_t(a^*a)$ with $t \in T \setminus S$ are invertible, whereas the operators $W_s (a^*a)$ with $s \in S$ are self-adjoint Fredholm operators with index zero. Hence, there are compact operators $K_s \in K(H_s)$ such that the operators $W_s(a^*a) + K_s$ are invertible. Let $j_s \in \cJ_s$ such that $W_s(j_s) = K_s$ and set $j := \sum_{s \in S} j_s$. Due to the separation property $(b)$ in Theorem \ref{tb1.6}, we then have
\[
W_t(a^*a + j) = \left\{
\begin{array}{ll}
W_t(a^*a) + K_t & \mbox{for} \; t \in S \\
W_t(a^*a)       & \mbox{for} \; t \in T \setminus S.
\end{array} \right.
\]
Thus, the operators $W_t(a^*a + j)$ are invertible for every $t \in T$, and their inverses are uniformly bounded by assumption. Since the family $(W_t)_{t \in T}$ is spectral, we get the invertibility of $a^*a + j$ in $\cA$. Further, $a^*a + j$ and $a^*a$ belong to the same coset modulo $\cJ$, whence the invertibility of $a + \cJ$ from the left-hand side. Repeating these arguments with $aa^*$ in place of $a^*a$ yields the invertibility of $a + \cJ$ from the right-hand side. The proof of assertion $(c)$ is similar to that of assertion $(a)$. \hfill \qed
\begin{coro} \label{cb1.14}
Let the notations and assumptions be as in Theorem $\ref{tb1.6}$, and suppose moreover that $(W_t)_{t \in T}$ is a spectral family for $\cA$. Put $W := \Pi_{t \in T} W_t$ and $H := \oplus_{ \in T} H_t$. Then a coset $a + \cJ$ is invertible if and only if $W(a)$ is a Fredholm operator on $L(H)$.
\end{coro}
Thus, under the hypotheses of the corollary, to study invertibility in $\cA/\cJ$ is the same as to study Fredholm theory on $L(H)$. \\[3mm]
{\bf Proof.} The spectral property of $(W_t)_{t \in T}$ is equivalent to the injectivity of $W$, considered as a mapping from $\cA$ to the product of the $C^*$-algebras $L(H_t)$. This product can be viewed as a $C^*$-subalgebra of the algebra $L(H)$, and the ideal $\cJ$ can then be identified with $W({\cA}) \cap K(H)$. Since
\[
(W(\cA) + K(H))/K(H) \cong W(\cA)/(W(\cA) \cap K(H)) \cong W(\cA)/W(\cJ),
\]
the assertion follows. \hfill \qed
\begin{coro} \label{cb1.14a}
Let $\cJ$ be a dual algebra constituted by a family $(\cJ_t)_{t \in T}$ of elementary ideals $\cJ_t$ with associated irreducible representations $W_t$. \\[1mm]
$(a)$ The algebra $\cJ$ is unital if and only if $T$ is finite and each algebra $\cJ_t$ has finite dimension. In this case, $(W_t)_{t \in T}$ is a strictly spectral family for $\cJ$. \\[1mm]
$(b)$ Let $\cJ$ be not unital and $\cA = \sC e + \cJ$ its minimal unitization. Then the family $(W_t)_{t \in T}$ is spectral for $\cA$, and this family is strictly spectral for $\cA$ if at least one of the algebras $\cJ_t$ has infinite dimension.
\end{coro}
Note that $\cJ$ cannot be unital if one of the algebras $\cJ_t$ has infinite dimension (otherwise the identity operator would be compact). \\[3mm]
{\bf Proof.} The proof of assertion $(a)$ is clear. Let $\cJ$ be non-unital, and let $a$ be an element of $\cA$ such that all $W_t(a)$ are invertible. The element $a$ can be uniquely written as $a = \gamma e + j$ with $\gamma \in \sC$ and $j \in \cJ$. Then $W_t(a) = \gamma I + W_t(j)$ is invertible for every $t \in T$. If $t \in T$ is such that $\cJ_t$ has infinite dimension, then $W_t(j)$ is a compact operator on an infinite-dimensional Hilbert space, which implies that $\gamma \neq 0$. But then $a$ is invertible modulo $\cJ$, and the Lifting theorem $\ref{tb1.6}$ implies the invertibility of $a$. Hence, $(W_t)_{t \in T}$ is a strictly spectral for $\cA$ in this case.

To prove the first assertion of $(b)$ we can thus assume that all elementary algebras $\cJ_t$ have finite dimension. Then $T$ must be infinite (otherwise $\cJ$ would be unital). Let now $a = \gamma e + j \in \cA$ be an element such that all $W_t(a)$ are invertible and the norms of their inverses are uniformly bounded, say by $M$. We claim that $\gamma \neq 0$ also in this case. Assume that $\gamma = 0$. Then
\[
1 \le \|W_t(\gamma e + j)\| \, \|W_t(a)^{-1}\| \le M \, \|W_t(j)\|
\]
for every $t \in T$. Since $\cJ$ is a dual algebra and $T$ is infinite, there is a $t \in T$ such that $\|W_t(j)\| < M$. For this $t$, the previous estimate yields a contradiction. \hfill \qed
\section{From dual ideals to Moore-Penrose ideals} \label{ssb1.4}
Here we prove the implication $(a) \Rightarrow (b)$ in Theorem \ref{tb3.90}, which we formulate as a separate result for further reference.
\begin{theo} \label{tb1.15}
Every proper dual ideal of a unital $C^*$-algebra $\cA$ is a Moore-Penrose ideal of $\cA$.
\end{theo}
{\bf Proof.} Let $\cJ$ be a proper dual ideal of $\cA$. We use the notations from the Lifting theorem \ref{tb1.6}. Thus, the ideal $\cJ$ is composed by its elementary ideals $\cJ_t$, and the homomorphisms $W_t : \cA \to L(H_t)$ lift the ideals $\cJ_t$.

Let $a \in \cA$ be invertible modulo $\cJ$. Then $W_t(a)$ is Fredholm for all $t \in T$ and $W_t(a)$ is invertible for all but finitely many $t$ by Theorem \ref{tb1.6} $(d)$. Let $P_t$ denote the orthogonal projection form $H_t$ onto the kernel of $W_t(a)$. These projections are compact, even of finite rank, and all but finitely many of them are zero. Thus, there are (uniquely determined) elements $p_t \in \cJ_t$ with $W_t(p_t) = P_t$. The sum $p := \sum_{t \in T} p_t$ is well defined, since only finitely many of the $p_t$ are non-zero, and $p \in \cJ$ is a projection. By the separation property,
\[
W_t(a^*a p) = W_t(a)^* W_t(a) W_t(p_t) = W_t(a)^* W_t(a) P_t = 0
\]
for every $t \in T$. Since $a^*a p \in \cJ$, and since $(W_t)_{t \in T}$ is a spectral family for the algebra $\sC e + \cJ$ by Corollary \ref{cb1.14a} $(b)$, we conclude from (\ref{e100}) that $a^*ap = 0$. Further, due to the choice of $P_t$, the operator
\[
W_t(a^*a + p) = W_t(a)^* W_t(a) + W_t(p_t) = W_t(a)^* W_t(a) + P_t
\]
is invertible for every $t \in T$. Since $a^*a + p$ is invertible modulo $\cJ$, the Lifting theorem \ref{ta6.8} implies the invertibility of $a^*a + p$. Thus, by Theorem \ref{t31.2}, $a$ is Moore-Penrose invertible in $\cA$, and $p$ is the associated Moore-Penrose projection. \hfill \qed \\[3mm]
It seems to be appropriate to give a more conceptual proof of Theorem \ref{tb1.15} which is based on some auxiliary facts which are of their own interest. The first fact is a modification of Lemma 6.3.1 in \cite{RSS2} (the $N$ ideals lemma) which holds for generalized invertibility.
\begin{lemma} \label{lb1.15a}
Let $\cJ_1, \, \cJ_2$ be ideals in a unital algebra $\cA$. If $a \in \cA$ is generalized invertible modulo $\cJ_1$ and one-sided invertible modulo $\cJ_2$, then $a$ is generalized invertible modulo $\cJ_1 \cap \cJ_2$.
\end{lemma}
{\bf Proof.} Let $e$ be the identity element of $\cA$, and let $b, \, c \in \cA$ be such that $aba - a =: j_1 \in \cJ_1$ and $e - ca =: j_2 \in \cJ_2$. Then
\[
j_1 j_2 = (aba-a)(e - ca) = a(b + c - bac)a - a \in \cJ_1 \cJ_2 \subseteq \cJ_1 \cap \cJ_2.
\]
Hence, $a$ is generalized invertible modulo $\cJ_1 \cap \cJ_2$. If $a$ is one-sided invertible modulo $\cJ_2$ from the other side, the proof is analogous. \hfill \qed \\[3mm]
The next fact settles the simplest case of Theorem \ref{tb1.15}, when $\cJ$ is elementary.
\begin{prop} \label{pb1.15b}
Every proper elementary ideal of a unital $C^*$-algebra $\cA$ is a Moore-Penrose ideal of $\cA$.
\end{prop}
{\bf Proof.} Let $\cJ$ be a proper elementary ideal of $\cA$ and $W$ a $^*$-isomorphism from $\cJ$ onto the ideal $K(H)$ of the compact operators on some Hilbert space $H$. Then $(W, \, H)$ is an irreducible representation of $\cJ$. This representation can be uniquely (up to unitary equivalence) extended to an irreducible representation of $\cA$ which we denote by $(W, \, H)$ again.

Let $a \in \cA$ be invertible modulo $\cJ$. Then $W(a)$ is invertible modulo $K(H)$, hence a Fredholm operator. In particular, $W(a)$ has a closed range and is, thus, Moore-Penrose invertible in $L(H)$. Then, by Theorem \ref{t31.2} $(f)$, $W(a)$ is Moore-Penrose invertible already in $W(\cA)$. Thus, $a$ is Moore-Penrose invertible in $\cA$ modulo $\ker W$. Since $a$ is also invertible modulo $\cJ$, it is generalized invertible modulo $\ker W \cap \cJ$ by the previous lemma. Moreover, $\ker W \cap \cJ = \{0\}$, since $W$ acts as an isometry on $\cJ$. Hence, $a$ is generalized invertible and , by Theorem \ref{t31.2}, also Moore-Penrose invertible in $\cA$. \hfill \qed \\[3mm]
The final auxiliary fact we need is that a direct sum of Moore-Penrose ideals is a Moore-Penrose ideal.
\begin{prop} \label{pb1.15c}
For every element $t$ of a non-empty set $T$, let $\cJ_t$ be a Moore-Penrose ideal of a unital $C^*$-algebra $\cA_t$. Then the direct sum $\oplus_{t \in T} \cJ_t$ is a Moore-Penrose ideal of the direct product $\Pi_{t \in T} \cA_t$.
\end{prop}
{\bf Proof.} Let $e_t$ denote the identity element of $\cA_t$. Let $a \in \Pi_{t \in T} \cA_t$ be invertible modulo $\oplus_{t \in T} \cJ_t$. Then there are elements $b \in \Pi_{t \in T} \cA_t$ and $j, \, k \in \oplus_{t \in T} \cJ_t$ such that $a_t b_t = e_t - j_t$ and $b_t a_t = e_t - k_t$ for every $t \in T$. Since $j, \, k \in \oplus_{t \in T} \cJ_t$, there is a finite subset $S$ of $T$ such that $\|j_t\| < 1/2$ and  $\|k_t\| < 1/2$ for every $t \in T \setminus S$. Thus, by Neumann series, $a_t$ is invertible in $\cA_t$ for $t \in T \setminus S$, and $\sup_{T \setminus S} \|a_t^{-1}\| < \infty$.

Let now $s \in S$. Then $a_s$ is invertible modulo the Moore-Penrose ideal $\cJ_s$; hence, $a_s$ is Moore-Penrose invertible in $\cA_s$. Since $S$ is finite,
\[
\max \, \{ \sup_{T \setminus S} \|a_t^{-1}\|, \, \max_{s \in S} \|a_s^\dagger\| \} < \infty.
\]
Thus, the function which sends $t \in T$ to $a_t^{-1}$ if $t \in T \setminus S$ and to $a_t^\dagger$ if $t \in S$ belongs to $\Pi_{t \in T} \cA_t$. This function is the Moore-Penrose inverse of $a$. \hfill \qed \\[3mm]
{\bf Alternate proof of Theorem \ref{tb1.15}.} Let $\cJ$ be a proper dual ideal of $\cA$ which is composed by its elementary ideals $\cJ_t$. We shall define a set $S$ and for every $s \in S$ a unital $C^*$-algebra $\cB_s$ with a dual ideal $\cI_s$ such that $\cI := \oplus_{s \in S} \cI_s$ is $^*$-isomorphic to $\cJ$ and $\cA$ is $^*$-isomorphic to a $C^*$-subalgebra of $\Pi_{s \in S} \cB_s$ which contains $\cI$.

Let $S := T \cup \{\ast\}$, define $\cB_s := \cA/\ann \cJ_s$ if $s \in T$ and $\cB_\ast := \cA/\cJ$, and let $\pi : \cA \to \Pi_{s \in S} \cB_s$ be the product of the canonical homomorphisms $\pi_s : \cA \to \cA/\ann \cJ_s$ if $s \in T$ and of the canonical homomorphism $\pi_\ast : \cA \to \cA/\cJ$. Since $\ker \pi = \ann \cJ \cap \cJ = \{0\}$, the mapping $\pi$ is an isometry (this is essentially the Lifting theorem \ref{ta6.8}). Thus, $\cA$ can be identified with the $C^*$-subalgebra $\pi(\cA)$ of $\Pi_{s \in S} \cB_s$. Moreover, $\pi$ sends the ideal $\cJ \cong \oplus_{s \in T} \cJ_t$ to
\[
\oplus_{t \in T} (\cJ + \ann \cJ_t)/\ann \cJ_t \oplus \cJ/\cJ \cong \oplus_{t \in T} \cJ/(\cJ \cap \ann \cJ_t) \oplus \{0\}.
\]
The ideals $\cI_t := \cJ/(\cJ \cap \ann \cJ_t)$ are $^*$-isomorphic to $\cJ_t$, hence elementary, and $\cJ$ is $^*$-isomorphic to $\oplus_{t \in T} \cI_t \oplus \cI_\ast$ with $\cI_\ast := \{0\}$ via the mapping $\pi$.

Now we can argue as follows. The ideals $\cI_s$ are Moore-Penrose ideals by Proposition \ref{pb1.15b} for $s \in T$, and $\cI_\ast$ is clearly also a Moore-Penrose ideal. Hence, by Proposition \ref{pb1.15c}, $\pi(\cJ) = \oplus_{s \in S} \cI_s$ is a Moore-Penrose ideal in $\Pi_{s \in S} \cB_s$ and, due to inverse closedness, also in $\pi(\cA)$. Since $\pi : \cA \to \pi(\cA)$ is a $^*$-isomorphism, $\cJ$ is a Moore-Penrose ideal in $\cA$. \hfill \qed
\section{From Moore-Penrose ideals to the compact spectral property} \label{ssb1.5}
Recall that an ideal $\cJ$ of a unital $C^*$-algebra has the compact spectral property if every self-adjoint element of $\cJ$ has 0 as the only possible accumulation point of its spectrum.
\begin{prop} \label{pb1.20a}
Every Moore-Penrose ideal in a unital $C^*$-algebra has the compact spectral property.
\end{prop}
{\bf Proof.} Let $a$ be a self-adjoint element of a Moore-Penrose ideal $\cJ$, and let $\lambda \in \sigma(a) \setminus \{0\}$. Then $a - \lambda e$ is invertible modulo $\cJ$ and thus, by the definition of a Moore-Penrose ideal, Moore-Penrose invertible. By Theorem \ref{t31.2}), 0 is an isolated point of the spectrum of $a - \lambda e$, hence, $\lambda$ is an isolated point of the spectrum of $a$. Thus, 0 is the only possible accumulation point in the spectrum of $a$. \hfill \qed
\section{From the compact spectral property to compact elements} \label{ssb1.6}
The aim of this section is to prove the implication $(c) \Rightarrow (e)$ in Theorem \ref{tb3.90}. 
\begin{theo} \label{tb1.20a1}
Let $\cJ$ be a closed ideal of a unital $C^*$-algebra $\cA$. If $\cJ$ has the compact spectral property, then $\cC(\cJ) = \cJ$. 
\end{theo}
We split the proof of Theorem \ref{tb1.20a1} into several steps which we formulate as separate assertions. The idea of the proof is that if a projection $q \in \cJ$ is not of finite algebraic rank, then one can perturb $q$ to obtain a self-adjoint element $a \in \cJ$ the spectrum of which contains a cloud of points with $1$ as accumulation point.
\begin{lemma} \label{pb2.2a}
Let $\cA$ be a $C^*$-algebra and $\cJ$ a closed ideal of $\cA$. If $k \in \cJ$ is an element of algebraic rank one in $\cJ$, then $k$ is also of algebraic rank one in $\cA$.
\end{lemma}
{\bf Proof.} Let $a \in A$ and $k \in \cJ$ an algebraic rank one element in $\cJ$. Let $(e_t)_{t \in T}$ be an approximate identity in $\cJ$. Then $k = \lim_{t \in T} k e_t$, whence
\begin{equation} \label{eb2.2b}
kak =  \lim_{t \in T} \, (k e_t) a k =  \lim_{t \in T} k (e_t a) k.
\end{equation}
Since $e_t a \in \cJ$, there are complex numbers $\mu_t$ such that $k (e_t a) k = \mu_t k$. From (\ref{eb2.2b})
we conclude that the limit $\lim_{t \in T} \mu_t k$ exists. This limit is necessarily of the form $\mu k$ with a complex number $\mu$. Hence, $kak = \mu k$. \hfill \qed 
\begin{lemma} \label{lb1.20a2}
Let $\cJ$ be a closed ideal of a unital $C^*$-algebra with compact spectral property and $q$ a non-zero projection in $\cJ$. Then every decreasing chain $q = q_0 > q_1 > q_2 > \ldots$ of projections (with proper relations $>$) is finite.
\end{lemma}
{\bf Proof.} Suppose there is a (countably) infinite chain $q = q_0 > q_1 > q_2 > \ldots$ of projections. Then $q_n = q_n q q_n$; hence $q_n \in \cJ$ for every $n$. Let $(\alpha_n)$ be an arbitrary real sequence in $l^1(\sN)$ and put
\[
a := q - \sum_{n = 1}^\infty \alpha_n (q_n - q_{n+1}).
\]
The sequence converges absolutely since $\|q_n - q_{n+1}\| = 1$, and it defines a self-adjoint element $a \in \cJ$. We consider the identity element $e$ of $\cA$, the element $a$ and the projections $q_n$ as elements of a commutative $C^*$-subalgebra of $\cA$ with maximal ideal space $X$. The relation $q_n > q_{n+1}$ implies that if $q_{n+1} (x) = 1$ for some $x \in X$ then $q_n(x) = 1$, and that if $q_n(x) = 0$ for some $x \in X$ then $q_{n+1} (x) = 0$. Moreover, since $q_n > q_{n+1}$ properly, there is an $x_n \in X$ such that $q_n(x_n) = 1$ and $q_{n+1}(x_n) = 0$. The equality
\[
a(x_i) = q(x_i) - \sum_{n = 1}^\infty \alpha_n (q_n (x_i) - q_{n+1} (x_i)) = q(x_i) - \alpha_i q_i(x_i) = 1 - \alpha_i
\]
shows that $1 - \alpha_i$ belongs to the spectrum of $a$ for every $i \in \sN$. In particular, if the $\alpha_i$ are pairwise distinct, then $1$ is an accumulation point of the spectrum of $a$, which contradicts the hypothesis that $0$ is the {\em only} accumulation point of the spectrum of the self-adjoint element $a$ of $\cJ$. \hfill \qed
\begin{coro} \label{cb1.20a3}
Let $\cJ$ be a closed ideal of a unital $C^*$-algebra with compact spectral property and $q$ a non-zero projection in $\cJ$. Then there is a minimal projection $r$ with $r < q$.
\end{coro}
{\bf Proof.} If $q$ is not minimal, there is a non-zero projection $q_1$ with $q > q_1$. If also $q_1$ is not minimal, there is a projection $q_2$ with $q_1 > q_2$. By Lemma \ref{lb1.20a2}, the chain $q = q_0 > q_1 > q_2 > \ldots$ is finite. Thus, one of the $q_i$ is minimal. \hfill \qed
\begin{coro} \label{cb1.20a4}
Let $\cJ$ be a closed ideal of a unital $C^*$-algebra with compact spectral property and $q$ a projection in $\cJ$. Then $q$ is a finite sum of minimal projections in $\cJ$.
\end{coro}
{\bf Proof.} If $q$ is already minimal, then there is nothing to prove. If $q$ is not minimal, then there exists a minimal projection $q_1 < q$ by Corollary \ref{cb1.20a3}. Then the element $q - q_1$ is a projection. If $q - q_1$ is minimal, then $q = q_1 + (q - q_1)$ is a sum of two minimal projections, and we are done. If $q - q_1$ is not minimal then there is a minimal projection $q_2 < q - q_1$ by Corollary \ref{cb1.20a3} again. If $q - q_1 - q_2$ is minimal, then $q = q_1 + q_2 + (q - q_1 - q_2)$ is a sum of three minimal projections. Otherwise, we proceed in this way. Since the chain $q > q - q_1 > q - q_1 - q_2 > \ldots$ is finite by Lemma \ref{lb1.20a2}, this process terminates. Thus, there is a non-trivial minimal projection $q - q_1 - \ldots - q_r$. Then $q = q_1 + \ldots + q_r + (q - q_1 - \ldots - q_r)$ is the sum of $r+1$ non-trivial minimal projections. \hfill \qed
\begin{prop} \label{pb1.20a5}
Let $\cJ$ be a closed ideal of a unital $C^*$-algebra with compact spectral property and $j \in \cJ \setminus \{0\}$. Then $j^*j$ is a limit of finite linear combinations of minimal projections in $\cJ$.
\end{prop}
{\bf Proof.} Since $\cJ$ has the compact spectral property, the spectrum of $j^*j$ is at most countable and has $0$ as only possible accumulation points. Let $\lambda_1 > \lambda_2 > \ldots$ denote the non-zero points in $\sigma(j^*j)$. Note that this chain is either finite or converges to $0$.

Since $\lambda_1$ is an isolated point of $\sigma (j^*j)$, the function $\widehat{p_1}$ which is $1$ at $\lambda_1$ and $0$ on $\sigma(j^*j) \setminus \{\lambda_1\}$ is continuous on $\sigma(j^*j)$. Since $\sigma(j^*j)$ is homeomorphic to the maximal ideal space of the smallest unital $C^*$-subalgebra of $\cJ$ which contains $j^*j$, the function $\widehat{p_1}$ defines an element $p_1$ of that algebra. This element is a non-zero projection, the so-called spectral projection of $j^*j$ associated with $\lambda_1$. Note that $p_1 \in \cJ$. By Corollary \ref{cb1.20a4}, $p_1$ is a finite sum of minimal projections in $\cJ$. Thus if $j^*j = \lambda_1 p_1$, then we are already done.

If $j^*j \neq \lambda p$, then we repeat these arguments, but now for the element $j^*j - \lambda_1 p_1$. Thus, we let $p_2$ be the spectral projection of $j^*j - \lambda_1 p_1$ associated with the largest point $\lambda_2$ is the spectrum of that element. If $j^*j = \lambda_1 p_1 + \lambda_2 p_2$, then we are done again. Otherwise we proceed in this way and consider the spectral projection $p_3$ of $j^*j - (\lambda_1 p_1 + \lambda_2 p_2)$ associated with $\lambda_3$.
Since
\[
\|j^*j - (\lambda_1 p_1 + \ldots + \lambda_n p_n)\| = \lambda_{n+1} \to 0 \quad \mbox{as} \; n \to \infty,
\]
and since each projection $p_i$ is a finite sum of minimal projections in $\cJ$ by Corollary \ref{cb1.20a4}, the assertion follows. \hfill \qed \\[3mm]
Since every element of $\cJ$ is a linear combination of four positive elements, it follows from Proposition \ref{pb1.20a5} that in fact {\em every} element of $\cJ$ is a limit of finite linear combinations of minimal projections in $\cJ$. Thus, the proof of Theorem \ref{tb1.20a1} will be complete once we have shown the following. Recall from Lemma \ref{pb2.2a} that every element of algebraic rank one in $\cJ$ is also of algebraic rank one in $\cA$.
\begin{prop} \label{pb1.20a6}
Let $\cJ$ be a closed ideal of a unital $C^*$-algebra with compact spectral property. Then every minimal projection in $\cJ$ is of algebraic rank one in $\cJ$.
\end{prop}
{\bf Proof.} Let $p$ be a minimal projection in $\cJ$ and $k \in \cJ$. We have to show that $pkp = \alpha p$ with a complex number $\alpha$. Since $k$ can be written as a linear combination of four non-negative elements, we can suppose without loss that $k = j^*j$ with $j \in \cJ$. Let $\cC$ stand for the smallest closed subalgebra of $\cJ$ which contains the elements $p$ and $pj^*jp$. This algebra is commutative and unital (with $p$ as its identity element), and its maximal ideal space is homeomorphic to $\sigma_\cC (pj^*jp)$. Note that $\sigma_\cC (pj^*jp) \subseteq \sigma_\cA (pj^*jp)$.

Assume that $\sigma_\cC (pj^*jp)$ contains at least two points. Then one of these points, $\lambda_0$ say, is not 0. This point isolated in $\sigma_\cA (pj^*jp)$ by the compact spectral property of $\cJ$. Hence, $\lambda_0$ is also isolated in $\sigma_\cC (pj^*jp)$. Then the function $\hat{q}$ which is $1$ at $\lambda_0$ and $0$ at all other points of $\sigma_\cC (pj^*jp)$ is continuous on $\sigma_\cC (pj^*jp)$. This function defines an element $q$, which is a projection in $\cC$.

Since $p$ is the identity element of $\cC$, one has $q \le p$. But $q$ is neither equal to $0$ (since $\hat{q}(\lambda_0) = 1$) nor equal to $p$ (since $\sigma_\cC (pj^*jp)$ contains at least two points and $\hat{p}$ is equal to one on all of $\sigma_\cC (pj^*jp)$). This contradicts the minimality of $p$. Hence, the maximal ideal space $\sigma_\cC (qb^*bq)$ of $\cC$ consists of exactly one point. Then every element of $\cC$ is a multiple of the identity element $p$. Hence, there is an $\alpha \in \sC$ such that $pj^*jp = \alpha p$. \hfill \qed
\section{From compact elements to dual ideals} \label{ssb1.6b}
In Theorem 3.3 in \cite{Roc2} we showed that every closed ideal generated by an element of algebraic rank one is elementary. Let $k, \, r$ be rank one elements and $\id (k)$, $\id (r)$ the elementary ideals generated by them. Since elementary algebras do not possess non-trivial ideals, it follows that either $\id (k) = \id (r)$ or $\id (k) \cap \id (r) = \{0\}$. From this property it easily follows that the smallest ideal which contains all of these elementary ideals id dual.  
\begin{theo} \label{tb2.23a1}
Let $\cA$ be a unital $C^*$-algebra. Then $\cC(\cA)$ is a dual algebra. 
\end{theo}
In particular, the equality $\cC(\cJ) = \cJ$ implies the duality of $\cJ$, which settles the implication $(e) \Rightarrow (a)$ in \ref{tb3.90}. 

To finish the proof of that theorem, we still have to show that $(a) \Rightarrow (d) \Rightarrow (c)$.

Let $\cJ$ be a dual ideal. Then $\cJ$ is $^*$-isomorphic to a direct sum $\oplus_{t \in T} K(H_t)$ of ideals of compact operators on certain Hilbert spaces $H_t$. We let $H$ denote the orthogonal sum of the Hilbert spaces $H_t$ with $t \in T$ and consider $H_t$ as a closed subspace of $H$. Let $P_t$ be the orthogonal projection from $H$ to $H_t$. Then $P_t K(H) P_t$ is $^*$-isomorphic to $K(H_t)$, and the family $(P_t K(H) P_t)_{t \in T}$ of $C^*$-subalgebras of $K(H)$ is separated in the sense that
\[
(P_s K(H) P_s)(P_t K(H) P_t) = \{0\} \qquad \mbox{whenever} \; s \neq t 
\]
since the Hilbert spaces $H_t$ are pairwise orthogonal. Thus, the smallest closed $C^*$-subalgebra of $K(H)$ which contains all algebras $P_t K(H) P_t \cong K(H_t)$ is dual, and this algebra is $C^*$-isomorphic to $\cJ$. This settles the implication $(a) \Rightarrow (d)$, and the implication $(d) \Rightarrow (c)$ is clear, since every compact operator on a Hilbert space has $0$ as only possible accumulation point in its spectrum. \hfill \qed 
\section{Lifting properties of Moore-Penrose ideals} \label{ssb1.7}
It is immediate from their definition that a Moore-Penrose ideal $\cJ$ of a unital $C^*$-algebra $\cA$ lifts Moore-Penrose invertible elements in the following sense: If a coset $a + \cJ$ is Moore-Penrose invertible in $\cA/\cJ$, then there is a $j \in \cJ$ such that $a + j$ is Moore-Penrose invertible in $\cA$. The following is then an immediate consequence of Theorem \ref{tb1.15}.
\begin{coro} \label{cb1.20}
Dual ideals of unital $C^*$-algebras lift Moore-Penrose invertible elements.
\end{coro}
Note in that connection that not every ideal that lifts Moore-Penrose invertible elements is a Moore-Penrose (= dual) ideal. For example, the ideal $\cJ := \{f \in C([0, \, 2]): f|_{[0, \, 1]} = 0\}$ of $C([0, \, 2])$ is certainly not dual, but it lifts Moore-Penrose invertible elements. Indeed, if $f + \cJ$ is Moore-Penrose invertible, then either $f$ is identically zero on $[0, \, 1]$, or it is non-zero on all of $[0, \, 1]$. In the first case, put $j := -f$. Then $j \in \cJ$ and $a + j = 0$ is a Moore-Penrose invertible lift of $f + \cJ$. In the second case, let $g$ be the continuation from $f|_{[0, \, 1]}$ to a function in $C([0, \, 2])$ which is constant on $[1, \, 2]$. Then $g-f \in \cJ$; hence, $g$ is an invertible lift of $f + \cJ$.

Another feature of dual ideals is that they lift projections in the following sense: every coset $a + \cJ$ which is a projection in $\cA/\cJ$ contains a representative which is a projection in $\cA$. Not every closed ideal of a $C^*$-algebra has the Moore-Penrose or the projection lifting property. For example, let $a \in C([0, \, 1])$ be a function with $a(0) = 0$ and $a(1) = 1$ and put $\cJ := \{f \in C([0, \, 1]) : f(0) = f(1) = 0 \}$. Then the coset $a + \cJ$ is a projection in $C([0, \, 1])/\cJ$, and $a + \cJ$ is Moore-Penrose invertible (in fact, every projection $q$ is Moore-Penrose invertible, and $q^\dagger = q$). But $a + \cJ$ can neither be lifted to a Moore-Penrose invertible function nor to a projection in $C([0, \, 1])$.
\begin{theo} \label{tb1.4}
Dual ideals of unital $C^*$-algebras lift projections.
\end{theo}
{\bf Proof.} Let $\cA$ be a $C^*$-algebra with identity element $e$ and $\cJ$ a dual ideal of $\cA$. Let $a \in \cA$ be such that $a + \cJ$ is a projection in $\cA/\cJ$. Since $a + \cJ$ is self-adjoint, we have $a - a^* \in \cJ$, hence $a - (a + a^*)/2 = (a - a^*)/2 \in \cJ$. Thus, the coset $a + \cJ$ contains a self-adjoint element. Let $a$ be this element.

Since $a + \cJ$ is an idempotent, the element $a - a^2 =: j$ is in $\cJ$. By Theorem \ref{tb3.90}, the point $0$ is the only possible accumulation point of the spectrum of $j$. Hence, the only possible accumulation points of the spectrum of $a$ are the points $0$ and $1$. In particular, there are only finitely many points $\lambda$ in $\sigma(a)$ with $|\lambda| \ge 1/2$ and $|1 - \lambda| \ge 1/2$, and these points are isolated in $\sigma(a)$. Let $\lambda$ be one of these points. Then the function
\[
r : \sigma(a) \to \sC, \quad r (x) := \left\{
\begin{array}{lll}
\lambda & \mbox{if} & x = \lambda \\
0 & \mbox{if} & x \in \sigma(a) \setminus \{ \lambda \}
\end{array} \right.
\]
is continuous and, thus, it defines an element of the $C^*$-subalgebra of $\cA$ generated by $a$ and $e$. We denote this element by $r$ again. Since $\lambda - \lambda^2 \neq 0$, we conclude from
\[
r (\lambda - \lambda^2) = r (a - a^2) = rj \in \cJ
\]
that $r \in \cJ$. Thus, the self-adjoint elements $a$ and $a - r$ lie in the same coset modulo $\cJ$, and $\sigma (a-r) = \sigma(a) \setminus \{ \lambda \}$. Repeating these arguments we obtain a self-adjoint element $b$ with $a-b \in \cJ$ such that
\[
\sigma(b) \subset \{ x \in \sR : |x| < 1/2 \} \cup \{ x \in \sR : |1-x| < 1/2 \} =: \Sigma_0 \cup \Sigma_1.
\]
Let $\cB$ denote the smallest closed subalgebra of $\cA$ which contains the elements $b$ and $e$. Then $\cB$ is a commutative $C^*$-algebra with maximal ideal space homeomorphic to $\sigma(b)$. The function $p : \sigma(b) \to \sR$ which is 0 on $\sigma(b) \cap \Sigma_0$ and 1 on $\sigma(b) \cap \Sigma_1$ is continuous on $\sigma(b)$. It defines a projection in $\cB$ which we denote by $p$ again. It remains to show that $b-p \in \cJ$. Let $c : \sigma(b) \to \sR$ be the function with $c(x) = 1/(1-x)$ if $x \in \sigma(b) \cap \Sigma_0$ and  $c(x) = -1/x$ if $x \in \sigma(b) \cap \Sigma_1$. This function is continuous on $\sigma(b)$, so it defines an element of $\cB$, and the equalities
\[
(b-p)(x) = \left\{
\begin{array}{ll}
x = \frac{1}{1-x} (x - x^2) = c(x) (b-b^2)(x) & \mbox{if} \; x \in \sigma(b) \cap \Sigma_0 \\
x - 1 = - \frac{1}{x} (x - x^2) = c(x) (b-b^2)(x) & \mbox{if} \; x \in \sigma(b) \cap \Sigma_1
\end{array} \right.
\]
show that $b-p = c(b-b^2)$. Since $b - b^2 \in \cJ$ we conclude that $b-p \in \cJ$, i.e., the projection $p$ is a representative of the coset $b + \cJ = a + \cJ$. \hfill \qed \\[3mm]
It is interesting to note that the projection lifting property of dual ideals follows already from their Moore-Penrose lifting property, which we observed in Corollary \ref{cb1.20}.
\begin{prop} \label{pb1.4a}
Let $\cJ$ be a closed ideal of a unital $C^*$-algebra. If $\cJ$ lifts Moore-Penrose invertible elements, then $\cJ$ lifts projections.
\end{prop}
{\bf Proof.} Let $p$ be an element of $\cA$ for which $p + \cJ$ is a projection. Then $p + \cJ$ is Moore-Penrose invertible. Since $\cJ$ lifts Moore-Penrose invertible elements, there is a Moore-Penrose invertible element $\tilde{p} \in \cA$ such that $p + \cJ = \tilde{p} + \cJ$. Let $r \in \cA$ be the Moore-Penrose projection of $\tilde{p}$, thus $\tilde{p}^* \tilde{p} + r$ is invertible in $\cA$ and $\tilde{p} r = 0$. Then
\[
(\tilde{p} + \cJ)^* (\tilde{p} + \cJ) + (r + \cJ) = (p + \cJ)^* (p + \cJ) + (r + \cJ) = p + r + \cJ
\]
is invertible in $\cA/\cJ$, and
\[
(\tilde{p} + \cJ) (r + \cJ) = (p + \cJ) (r + \cJ) = pr + \cJ = 0 + \cJ.
\]
The last identity implies that $p + r + \cJ$ is a projection in $\cA/\cJ$. As we have seen before, this projection is invertible. Hence, $p + r + \cJ = e + \cJ$ or, equivalently, $p + \cJ = (e-r) + \cJ$. Thus, $e-r$ is a projection which lifts $p + \cJ$. \hfill \qed \\[3mm]
The following example shows that the converse of Proposition \ref{pb1.4a} is wrong.
\begin{example} \label{exb1.4b}
Let $\overline{\sD} := \{ z \in \sC : |z| \le 1\}$ and $\sT := \{ z \in \sC : |z| = 1\}$, and put $\cA := C(\overline{\sD})$ and $\cJ = \{f \in C(\overline{\sD}) : f|_\sT = 0 \}$. The ideal $\cJ$ lifts projections. Indeed, let $p$ be a function in $C(\overline{\sD})$ such that $p + \cJ$ is a projection. Then $p(x) - p(x)^2 = 0$ on $\sT$, whence $p(x) \in \{0, \, 1\}$ for $x \in \sT$. Since $\sT$ is connected and $p$ is continuous, this implies that either $p(x) = 0$ for all $t \in \sT$ or $p(x) = 1$ for all $t \in \sT$. In the first case, the function $z \mapsto 0$ on $\overline{\sD}$ is a projection in $\cA$ which lifts $p + \cJ$, in the second case the function $z \mapsto 1$ does the job.

We will now see a function $a$ for which $a + \cJ$ is Moore-Penrose invertible (in fact, invertible), but the coset $a + \cJ$ contains no Moore-Penrose invertible representative. Let $a : z \mapsto z$. Then $\overline{a} (x) a(x) - 1 = 0$ for $x \in \sT$. Hence, $\overline{a} a + \cJ = 1 + \cJ$, and the coset $a + \cJ$ is invertible modulo $\cJ$. We show that there is no Moore-Penrose invertible function $\tilde{a} \in \cA$ with $a + \cJ = \tilde{a} + \cJ$. Suppose there is such a function. Then, by Theorem \ref{t31.2}, either the function $\tilde{a}$ is invertible, or $0$ is an isolated point of the spectrum of $\tilde{a} a$. The latter is impossible since the spectrum of $\tilde{a} a$ is the connected set $(\tilde{a} a) (\overline{\sD})$. Thus, if $0$ is a isolated point of this set, then $(\tilde{a} a) (\overline{\sD}) = \{0\}$. This implies that $\tilde{a}$ is the zero function, which is not a representative of $a + \cJ$. Hence, $\tilde{a}$ must be invertible on $\overline{\sD}$.

For $r \in [0, \, 1]$ and $x \in \sT$, put $a_r(x) = \tilde{a} (rx)$. Since $\tilde{a}$ is invertible on $\overline{\sD}$, the functions $a_r$ do not vanish and have, thus, a well defined winding number. Since the winding number depends continuously on $r$, it is independent of $r$. But the winding numbers of $a_0$ and $a_1$ are $0$ and $1$, respectively. This contradiction shows that $\cJ$ does not lift Moore-Penrose invertible elements. \hfill \qed \end{example}
Applied to the dual ideal $K(H)$ of $L(H)$, the following result states that every operator which is Moore-Penrose invertible modulo $K(H)$ is the sum of an operator with closed range and a compact operator.
\begin{theo} \label{tb1.16}
Let $\cA$ be a unital $C^*$-algebra, $\cJ$ a dual ideal of $\cA$, and $a \in \cA$. The coset $a + \cJ$ is Moore-Penrose invertible in $\cA/\cJ$ if and only if $a$ is the sum of an Moore-Penrose invertible element and an element in $\cJ$.
\end{theo}
{\bf Proof.} If $a \in \cA$ is the sum of a Moore-Penrose invertible element and an element in $\cJ$ then the coset $a + \cJ$ is clearly Moore-Penrose invertible in $\cA/\cJ$.

Conversely, let $a + \cJ$ be Moore-Penrose invertible in $\cA/\cJ$. Then $a^*a + \cJ$ is Moore-Penrose invertible in $\cA/\cJ$ by Theorem \ref{t31.2}. Let $\cB$ stand for the smallest closed subalgebra of $\cA$ which contains $a^*a$ and the identity element $e$ of $\cA$. Then the coset $a^*a + \cJ$ is Moore-Penrose invertible in $(\cB + \cJ)/\cJ$ by Theorem \ref{t31.2} $(f)$, whence the Moore-Penrose invertibility of the coset $a^*a + (\cB \cap \cJ)$ in $\cB/(\cB \cap \cJ)$ (the algebras $(\cB + \cJ)/\cJ$ and $\cB/(\cB \cap \cJ)$ are canonically isomorphic). Let $\pi \in \cB$ be an element such that $\pi + (\cB \cap \cJ)$ is the Moore-Penrose projection of $a^*a + (\cB \cap \cJ)$. Thus, $a^*a + \pi$ is invertible modulo $\cB \cap \cJ$, and $a^*a \pi \in \cB \cap \cJ$.

Being a subalgebra of a dual ideal, the algebra $\cA_0 \cap \cJ$ is dual by Theorem \ref{tb3.90} $(d)$, and dual ideals lift projections by Theorem \ref{tb1.4}. Hence, $\pi = p + r$ where $p$ is a projection in $\cB$ and $r \in \cB \cap \cJ$. Since $a^*a + p$ is invertible modulo $\cB \cap \cJ$ and $a^*a p \in \cB \cap \cJ$, there are elements $b \in \cB$ and $k, \, l \in \cB \cap \cJ$ such that
\begin{equation} \label{eb1.17}
b(a^*a + p) = (a^*a + p) b = e + k \quad \mbox{and} \quad a^*ap = l.
\end{equation}
Set $q := e - p$. Multiplying the first identity in (\ref{eb1.17}) by $q$ from both sides and employing the fact that all elements under consideration belong to the commutative algebra $\cB$, one easily gets
\begin{equation} \label{eb1.18}
(qbq) \, (qa^*aq) = (qa^*aq) \, (qbq) = q + qkq.
\end{equation}
Since $k \in \cJ$, this implies that $qa^*aq + p$ is invertible modulo in $\cB/(\cB \cap \cJ)$.

Now we move into the context of the Lifting theorem \ref{tb1.6}, which we apply to the dual ideal $\cJ$. With the notation as in the lifting theorem, let $t \in T$. Applying the homomorphism $W_t$ to both sides of (\ref{eb1.18}) yields
\begin{equation} \label{eb1.19}
W_t (qbq) \, W_t (qa^*aq) = W_t(q) + W_t (qkq)
\end{equation}
with a compact operator $W_t (qkq)$ is compact. Let $S$ be the set of all $t \in T$ with $\|W_t (qkq)\| < 1/2$. For $t \in S$, the right hand side of (\ref{eb1.19}) is invertible by Neumann series, and the norm of its inverse is not greater than 2. Hence,
\[
(W_t(q) + W_t (qkq))^{-1} W_t (qbq) \, W_t (qa^*) \, W_t(aq) = W_t(q)
\]
with $\|(W_t(q) + W_t (qkq))^{-1} W_t (qbq) W_t(qa^*)\| \le 2 \, \|qbq a^*\|$. Since all occurring operators belong to the commutative algebra $W_t(\cB)$, one easily concludes that the operators $W_t (aq)$ are normally solvable and that the norms of their Moore-Penrose inverses are uniformly bounded with respect to $t \in S$.

Let now $t \in T \setminus S$. Then (\ref{eb1.19}) implies that $W_t (qa^*aq)$ is a Fredholm operator. From Theorem \ref{t31.2} we conclude that $W_t (aq)$ is normally solvable in this case, too. Since $T \setminus S$ is a finite set, one gets that $W_t (aq)$ is normally solvable for all $t \in T$ and that the norms of their Moore-Penrose inverses are uniformly bounded. In terms of the product homomorphism $W := \Pi_{t \in T} W_t : \cA \to \Pi_{t \in T} L(H_t)$, this amounts to saying that $W(aq)$ is Moore-Penrose invertible in the $\Pi_{t \in T} L(H_t)$, hence in $W(\cA)$. Thus, the coset $aq + \ker W$ is Moore-Penrose invertible in the quotient algebra $\cA/\ker W$.

The second equality in (\ref{eb1.17}) implies via the $C^*$-axiom that $ap \in \cJ$. Thus, since $a + \cJ$ is
Moore-Penrose invertible, the coset $aq + \cJ = a - ap + \cJ$ is Moore-Penrose invertible, too. Summarizing, we have found that $aq$ is Moore-Penrose invertible modulo the ideals $\ker W$ and $\cJ$. Then the pair $(aq + \ker W, \, aq + \cJ)$ is Moore-Penrose invertible in $(\cA/\ker W) \times (\cA/\cJ)$. The mapping
\[
V : \cA \to (\cA/\ker W) \times (\cA/\cJ), \quad a \mapsto (a + \ker W, \, a + \cJ)
\]
is a $C^*$-homomorphism. Hence, $V(\cA)$ is a $C^*$-subalgebra of $(\cA/\ker W) \times (\cA/\cJ)$, and $V(aq)$ is Moore-Penrose invertible in $V(\cA)$ by Theorem \ref{t31.2} $(f)$. From Corollary 6.3.3 in \cite{RSS2} we recall that $\ker W \cap \cJ = \{0\}$. Hence, $V$ is a $^*$-isometry, and the Moore-Penrose invertibility of $V(aq)$ in $V(\cA)$ implies that of $aq$ in $\cA$. Thus, $a = ap + aq$ is the sum of an element in $\cJ$ and a Moore-Penrose invertible element. \hfill \qed
{\small Author's address: \\[3mm]
Steffen Roch, Technische Universit\"at Darmstadt, Fachbereich Mathematik, Schlossgartenstrasse 7, 64289 Darmstadt, Germany. \\
E-mail: roch@mathematik.tu-darmstadt.de}

\begin{thebibliography}{11}
%
%
\bibitem{Atk1}
{\sc F. W. Atkinson}, Normal solvability of linear equations in normed spaces. -- Mat. Sbornik {\bf 28}(1951), 1, 3 -- 14 (Russian).
\bibitem{Arv1}
{\sc W. Arveson}, An Invitation to $C^*$-Algebras. -- Springer Verlag, New York - Heidelberg - Berlin 1976.
\bibitem{Ber1}
{\sc M. C. F. Berglund}, Ideal $C^*$-algebras. -- Duke Math. J. {\bf 40}(1973), 241 -- 257.
\bibitem{MiP1}
{\sc S. G. Michlin, S. Pr\"ossdorf}, Singul\"are Integraloperatoren. -- Aka\-de\-mie-Verlag, Berlin 1980 (Extended English transl.: Singular Integral Operators. -- Akademie-Verlag, Berlin 1986, and Springer-Verlag, Heidelberg 1986).
\bibitem{Roc2}
{\sc S. Roch}, Algebras of approximation sequences: Fredholmness. -- J. Oper. Theory {\bf 48}(2002), 121 -- 149.
\bibitem{RSS2}
{\sc S. Roch, P. A. Santos, B. Silbermann}, Non-commutative Gelfand Theories. A Tool-kit for Operator Theorists and Numerical Analysts. -- Universitext, Springer, London 2011.
%
%
\end{thebibliography}
\end{document}